\documentclass{gtart}

\def\ifplaintex{\expandafter\ifx\csname documentclass\endcsname\relax}

\def\gtp{{\mathsurround=0pt\it $\cal G\mskip-2mu$eometry \&\ 
$\cal T\!\!$opology $\cal P\!$ublications}}  

\def\recd{{\small Received:\qua\receiveddate\ifx\reviseddate\relax
\else\qquad Revised:\qua\reviseddate\fi\par}} 

\def\lognumber#1{\def\thelognumber{#1}}
\def\volumenumber#1{\def\thevolumenumber{#1}}
\def\volumeyear#1{\def\thevolumeyear{#1}}
\def\papernumber#1{\def\thepapernumber{#1}}
\def\pagenumbers#1#2{\def\startpage{#1}\def\finishpage{#2}}
\def\published#1{\def\publishdate{#1}}

\def\received#1{\def\receiveddate{#1}}

\def\accepted#1{\def\accepteddate{#1}}

\def\asciiaddress#1{\def\theasciiaddress{#1}}
\def\asciiemail#1{\def\theasciiemail{#1}}

\long\def\asciiabstract#1{\long\def\theasciiabstract{#1}}

\let\\\par\let\thelognumber\relax\let\thevolumenumber\relax
\let\thepapernumber\relax\let\thevolumeyear\relax\let\startpage\relax
\let\finishpage\relax\let\publishdate\relax\let\receiveddate\relax
\let\reviseddate\relax\let\accepteddate\relax\let\theasciititle\relax
\let\theasciiauthors\relax\let\theasciiaddress\relax
\let\theasciiabstract\relax

\let\theasciiemail\relax

\ifplaintex
\font\logobig=cmssbx10 scaled 3836
\font\logomed=cmssbx10 scaled 2557
\else
\font\logobig=cmssbx10 scaled 4200
\font\logomed=cmssbx10 scaled 2800
\fi

\long\def\makeagttitle{   
\count0=\startpage
\agt\hfill      
\hbox to 45truept{\vbox to 0pt{\vglue -13truept{\logomed A\kern -.37em{\logobig 
T}\kern -.38em G}\vss}\hss}
\break
{\small Volume \thevolumenumber\ (\thevolumeyear)
\startpage--\finishpage\nl
Published: \publishdate}

\vglue .25truein

{\parskip=0pt\leftskip 0pt plus
1fil\def\\{\par\smallskip}{\Large\bf\thetitle}\par\medskip} \vglue
0.05truein

{\parskip=0pt\leftskip 0pt plus 1fil\def\\{\par}{\sc\theauthors}
\par\medskip}%
 
\vglue 0.03truein 

{\small\leftskip 25truept\rightskip 25truept{\bf Abstract}\stdspace\theabstract

{\bf AMS Classification}\stdspace\theprimaryclass
\ifx\thesecondaryclass\relax\else; \thesecondaryclass\fi\par
{\bf Keywords}\stdspace \thekeywords\par}\vglue 7truept

}   

\ifplaintex
\font\phead=cmsl9 scaled 950
\font\pnum=cmbx10 scaled 913
\font\pfoot=cmsl9 scaled 950
\headline{\vbox to 0pt{\vskip -4.5mm\line{\small\phead\ifnum
\count0=\startpage ISSN 1472-2739 (on-line) 1472-2747 (printed)
\hfill {\pnum\folio}\else\ifodd\count0\def\\{ }%
\ifx\theshorttitle\relax\thetitle\else\theshorttitle\fi\hfill{\pnum\folio}
\else\def\\{ and }{\pnum\folio}\hfill\ifx\theshortauthors\relax\theauthors
\else\theshortauthors\fi\fi\fi}\vss}}
\footline{\vbox to 0pt{\vglue 0mm\line{\small\pfoot\ifnum\count0=\startpage
\copyright\ \gtp\hfill\else
\agt, Volume \thevolumenumber\ (\thevolumeyear)\hfill\fi}\vss}}
\else
\font=cmsl9 scaled 1050
\font\lnum=cmbx10 
\font=cmsl9 scaled 1050
\makeatletter
\def\@oddhead{{\small\lhead\ifnum\count0=\startpage ISSN 1472-2739 
(on-line) 1472-2747 (printed)\hfill {\lnum\number\count0}\else\ifodd\count0
\def\\{ }\ifx\theshorttitle\relax \thetitle \else\theshorttitle\fi\hfill
{\lnum\number\count0}\else\def\\{ and }{\lnum\number\count0}
\hfill\ifx\theshortauthors\relax 
\theauthors\else\theshortauthors\fi\fi\fi}}\def\@evenhead{\@oddhead}
\def\@oddfoot{\small\lfoot\ifnum\count0=\startpage\copyright\ \gtp\hfill\else
\agt, Volume \thevolumenumber\ (\thevolumeyear)\hfill\fi}
\def\@evenfoot{\@oddfoot}
\makeatother
\fi
\let\maketitlepage\makeagttitle
\let\makeshorttitle\maketitlepage
\let\maketitle\maketitlepage

\newwrite\gtoutfile
\long\gdef\makeheadfile{  
{\def\\{, }\def\s{ }
\immediate\openout\gtoutfile head.xxx
\immediate\write\gtoutfile{To: math@arxiv.org}
\immediate\write\gtoutfile{Subject: put OR rep NNNNN:ppppp}
\immediate\write\gtoutfile{--text follows this line--}
\immediate\write\gtoutfile{Proxy-for: \ifx\theasciiauthors\relax
\theauthors\else\theasciiauthors\fi\s<\ifx\theasciiemail\relax\theemail\else\theasciiemail\fi>}
\immediate\write\gtoutfile{\noexpand\\}
\immediate\write\gtoutfile{Authors: \ifx\theasciiauthors\relax
\theauthors\else\theasciiauthors\fi}
{\def\\{ }\immediate\write\gtoutfile{Title: \ifx\theasciititle\relax
\thetitle\else\theasciititle\fi}}
\immediate\write\gtoutfile{Subj-class: GT or SG, GR etc}
\immediate\write\gtoutfile{MSC-class: \theprimaryclass\ifx\thesecondaryclass\relax\else, \thesecondaryclass\fi}
\immediate\write\gtoutfile{Journal-ref: Algebr. Geom. Topol. \thevolumenumber\s
(\thevolumeyear) \startpage-\finishpage}
\immediate\write\gtoutfile{Comments: Published by Algebraic and
Geometric Topology at}
\immediate\write\gtoutfile{\s\s\s  http://www.maths.warwick.ac.uk/agt/AGTVol\thevolumenumber/agt-\thevolumenumber-\thepapernumber.abs.html}
\immediate\write\gtoutfile{\noexpand\\}
\immediate\write\gtoutfile{}
\ifx\theasciiabstract\relax
\immediate\write\gtoutfile{\theabstract}\else
\immediate\write\gtoutfile{\theasciiabstract}\fi
\immediate\write\gtoutfile{}
\immediate\write\gtoutfile{\noexpand\\}
\immediate\write\gtoutfile{}
\immediate\closeout\gtoutfile}}  

\def\maketitlepage{\makeagttitle\makeheadfile}
\let\makeshorttitle\maketitlepage
\let\maketitle\maketitlepage

\def\ifplaintex{\expandafter\ifx\csname documentclass\endcsname\relax}

\def\gtp{{\mathsurround=0pt\it $\cal G\mskip-2mu$eometry \&\ 
$\cal T\!\!$opology $\cal P\!$ublications}}  

\def\recd{{\small Received:\qua\receiveddate\ifx\reviseddate\relax
\else\qquad Revised:\qua\reviseddate\fi\par}} 

\def\lognumber#1{\def\thelognumber{#1}}
\def\volumenumber#1{\def\thevolumenumber{#1}}
\def\volumeyear#1{\def\thevolumeyear{#1}}
\def\papernumber#1{\def\thepapernumber{#1}}
\def\pagenumbers#1#2{\def\startpage{#1}\def\finishpage{#2}}
\def\published#1{\def\publishdate{#1}}

\def\received#1{\def\receiveddate{#1}}

\def\accepted#1{\def\accepteddate{#1}}

\def\asciiaddress#1{\def\theasciiaddress{#1}}
\def\asciiemail#1{\def\theasciiemail{#1}}

\long\def\asciiabstract#1{\long\def\theasciiabstract{#1}}

\let\\\par\let\thelognumber\relax\let\thevolumenumber\relax
\let\thepapernumber\relax\let\thevolumeyear\relax\let\startpage\relax
\let\finishpage\relax\let\publishdate\relax\let\receiveddate\relax
\let\reviseddate\relax\let\accepteddate\relax\let\theasciititle\relax
\let\theasciiauthors\relax\let\theasciiaddress\relax
\let\theasciiabstract\relax

\let\theasciiemail\relax

\ifplaintex
\font\logobig=cmssbx10 scaled 3836
\font\logomed=cmssbx10 scaled 2557
\else
\font\logobig=cmssbx10 scaled 4200
\font\logomed=cmssbx10 scaled 2800
\fi

\long\def\makeagttitle{   
\count0=\startpage
\agt\hfill      
\hbox to 45truept{\vbox to 0pt{\vglue -13truept{\logomed A\kern -.37em{\logobig 
T}\kern -.38em G}\vss}\hss}
\break
{\small Volume \thevolumenumber\ (\thevolumeyear)
\startpage--\finishpage\nl
Published: \publishdate}

\vglue .25truein

{\parskip=0pt\leftskip 0pt plus
1fil\def\\{\par\smallskip}{\Large\bf\thetitle}\par\medskip} \vglue
0.05truein

{\parskip=0pt\leftskip 0pt plus 1fil\def\\{\par}{\sc\theauthors}
\par\medskip}%
 
\vglue 0.03truein 

{\small\leftskip 25truept\rightskip 25truept{\bf Abstract}\stdspace\theabstract

{\bf AMS Classification}\stdspace\theprimaryclass
\ifx\thesecondaryclass\relax\else; \thesecondaryclass\fi\par
{\bf Keywords}\stdspace \thekeywords\par}\vglue 7truept

}   

\ifplaintex
\font\phead=cmsl9 scaled 950
\font\pnum=cmbx10 scaled 913
\font\pfoot=cmsl9 scaled 950
\headline{\vbox to 0pt{\vskip -4.5mm\line{\small\phead\ifnum
\count0=\startpage ISSN 1472-2739 (on-line) 1472-2747 (printed)
\hfill {\pnum\folio}\else\ifodd\count0\def\\{ }%
\ifx\theshorttitle\relax\thetitle\else\theshorttitle\fi\hfill{\pnum\folio}
\else\def\\{ and }{\pnum\folio}\hfill\ifx\theshortauthors\relax\theauthors
\else\theshortauthors\fi\fi\fi}\vss}}
\footline{\vbox to 0pt{\vglue 0mm\line{\small\pfoot\ifnum\count0=\startpage
\copyright\ \gtp\hfill\else
\agt, Volume \thevolumenumber\ (\thevolumeyear)\hfill\fi}\vss}}
\else
\font=cmsl9 scaled 1050
\font\lnum=cmbx10 
\font=cmsl9 scaled 1050
\makeatletter
\def\@oddhead{{\small\lhead\ifnum\count0=\startpage ISSN 1472-2739 
(on-line) 1472-2747 (printed)\hfill {\lnum\number\count0}\else\ifodd\count0
\def\\{ }\ifx\theshorttitle\relax \thetitle \else\theshorttitle\fi\hfill
{\lnum\number\count0}\else\def\\{ and }{\lnum\number\count0}
\hfill\ifx\theshortauthors\relax 
\theauthors\else\theshortauthors\fi\fi\fi}}\def\@evenhead{\@oddhead}
\def\@oddfoot{\small\lfoot\ifnum\count0=\startpage\copyright\ \gtp\hfill\else
\agt, Volume \thevolumenumber\ (\thevolumeyear)\hfill\fi}
\def\@evenfoot{\@oddfoot}
\makeatother
\fi
\let\maketitlepage\makeagttitle
\let\makeshorttitle\maketitlepage
\let\maketitle\maketitlepage

\newwrite\gtoutfile
\long\gdef\makeheadfile{  
{\def\\{, }\def\s{ }
\immediate\openout\gtoutfile head.xxx
\immediate\write\gtoutfile{To: math@arxiv.org}
\immediate\write\gtoutfile{Subject: put OR rep NNNNN:ppppp}
\immediate\write\gtoutfile{--text follows this line--}
\immediate\write\gtoutfile{Proxy-for: \ifx\theasciiauthors\relax
\theauthors\else\theasciiauthors\fi\s<\ifx\theasciiemail\relax\theemail\else\theasciiemail\fi>}
\immediate\write\gtoutfile{\noexpand\\}
\immediate\write\gtoutfile{Authors: \ifx\theasciiauthors\relax
\theauthors\else\theasciiauthors\fi}
{\def\\{ }\immediate\write\gtoutfile{Title: \ifx\theasciititle\relax
\thetitle\else\theasciititle\fi}}
\immediate\write\gtoutfile{Subj-class: GT or SG, GR etc}
\immediate\write\gtoutfile{MSC-class: \theprimaryclass\ifx\thesecondaryclass\relax\else, \thesecondaryclass\fi}
\immediate\write\gtoutfile{Journal-ref: Algebr. Geom. Topol. \thevolumenumber\s
(\thevolumeyear) \startpage-\finishpage}
\immediate\write\gtoutfile{Comments: Published by Algebraic and
Geometric Topology at}
\immediate\write\gtoutfile{\s\s\s  http://www.maths.warwick.ac.uk/agt/AGTVol\thevolumenumber/agt-\thevolumenumber-\thepapernumber.abs.html}
\immediate\write\gtoutfile{\noexpand\\}
\immediate\write\gtoutfile{}
\ifx\theasciiabstract\relax
\immediate\write\gtoutfile{\theabstract}\else
\immediate\write\gtoutfile{\theasciiabstract}\fi
\immediate\write\gtoutfile{}
\immediate\write\gtoutfile{\noexpand\\}
\immediate\write\gtoutfile{}
\immediate\closeout\gtoutfile}}  

\def\maketitlepage{\makeagttitle\makeheadfile}
\let\makeshorttitle\maketitlepage
\let\maketitle\maketitlepage

\lognumber{18}
\volumenumber{1}
\volumeyear{2001}
\papernumber{18}
\pagenumbers{369}{380}
\received{4 April 2001}
\accepted{1 June 2001}
\published{3 June 2001}

\usepackage{amssymb}
\usepackage{amsmath}

\theoremstyle{plain}
\newtheorem{theorem}{Theorem}[section]
\newtheorem*{thm}{Theorem}
\newtheorem*{quotethm}{Theorem}
\newtheorem{lemma}[theorem]{Lemma}
\newtheorem{corollary}[theorem]{Corollary}

\theoremstyle{definition}
\newtheorem{definition}{Definition}[section]

\theoremstyle{remark}
\newtheorem*{remark}{Remark}
\newtheorem*{acknowledgements}{Acknowledgements}

\newcommand{\del}{\partial}

\newcommand{\Z}{{\mathbb{Z}}}

\newcommand{\R}{{\mathbb{R}}}

\newcommand{\sm}{-}

\begin{document} 

\title{Topological geodesics and virtual rigidity}

\authors{Louis Funar\\Siddhartha Gadgil}

\address{Institut Fourier BP74, UMR 5582, 
Universite de Grenoble I\\38402 Saint-Martin-d'Heres cedex, France} 
\secondaddress{Department of Mathematics, 
SUNY at Stony Brook\\Stony Brook, NY 11794, USA}

\asciiaddress{Institut Fourier BP74, UMR 5582\\ 
Universite de Grenoble I\\38402 Saint-Martin-d'Heres cedex, France\\ 
Department of Mathematics\\ 
SUNY at Stony Brook,Stony Brook, NY 11794, USA}

\email{funar@fourier.ujf-grenoble.fr}
\secondemail{gadgil@math.sunysb.edu}
\asciiemail{funar@fourier.ujf-grenoble.fr, gadgil@math.sunysb.edu}

\begin{abstract}
We introduce the notion of a topological geodesic in a $3$-manif\-old.
Under suitable hypotheses on the fundamental group, for instance
word-hyperbolicity, topological geodesics are shown to have the useful
properties of, and play the same role in several applications as,
geodesics in negatively curved spaces.  This permits us to obtain
virtual rigidity results for $3$-manifolds.
\end{abstract}
\asciiabstract{
We introduce the notion of a topological geodesic in a 3-manifold.
Under suitable hypotheses on the fundamental group, for instance
word-hyperbolicity, topological geodesics are shown to have the useful
properties of, and play the same role in several applications as,
geodesics in negatively curved spaces.  This permits us to obtain
virtual rigidity results for 3-manifolds.}

\maketitle 

Geodesics in Riemannian manifolds with metrics of negative sectional
curvature play an essential role in geometry. We show here that, in
the case of $3$-dimensional manifolds, many crucial properties of
geodesics follow from a purely topological characterization in terms
of \emph{knotting}. In particular, we prove two results concerning the
\emph{virtual rigidity} of $3$-manifolds following the methods of
Gabai~\cite{Ga1}.

We introduce the notion of a \emph{topological geodesic} in a
$3$-manifold. We shall prove basic existence and uniqueness results
for topological geodesics under suitable hypotheses on the fundamental group.

Suppose henceforth that $M$ is a closed $3$-manifold with
word-hyperbolic (or semi-hyperbolic) $\pi_1(M)$.
 We refer to the next section for the 
definition of the semi-hyperbolicity, following Alonso and Bridson.
 We adopt the
convention here that finite groups are not semi-hyperbolic, hence all
closed $3$-manifolds we consider have infinite fundamental group
unless the opposite is explicitly stated.  In particular, the
universal cover $\widetilde{M}$ of $M$ is homeomorphic to $\R^3$ (see
\cite{BM,M}).

\begin{definition} 
An embedded curve $\gamma$ in $M$ is a topological geodesic if a
component $\widetilde\gamma$ of its inverse image in $\widetilde M$ is unknotted.
\end{definition}

\begin{remark} 
This is equivalent to saying that every component of its inverse image
is unknotted. One observes that such a component is either a circle 
when $\gamma$ is torsion in $\pi_1(M)$ or a proper line in $\R^3$. 
\end{remark}

Under the hypothesis that $M$ is word-hyperbolic (or semi-hyperbolic)
we have the following existence theorem.

\begin{thm}[See theorem~\ref{T:exist} and corollary~\ref{T:exist2}]
Let $M$ be an irreducible $3$-manifold with $\pi_1(M)$ word-hyperbolic
(or semi-hyperbolic). Then every conjugacy class in
$\pi_1(M)$ is represented by a topological geodesic.
\end{thm}

If one further assumes that $\pi_1(M)$ is residually finite, we have an
uniqueness result.

\begin{thm}[See theorem~\ref{T:unique}]
Let $M$ be an irreducible $3$-manifold with $\pi_1(M)$
word-hyperbolic and residually finite. Suppose $c$ and $c'$
are homotopic topological geodesics in $M$ representing a primitive
class in $\pi_1(M)$ (i.e., not a multiple of any other class), then
there exists a finite cover $M'$ of $M$ such that $c$ and $c'$ lift to
isotopic curves in $M'$.
\end{thm}

In the case of a geodesic $\gamma$ in a Riemannian manifold, the
exponential map is a surjection on a neighbourhood of $\gamma$,
allowing one to construct a tubular neighbourhood. In the case of
negative sectional curvature one can do more. Namely, a
Hadamard-Cartan type argument allows one to construct thick tubular
neighbourhoods when the injectivity radius is large. Thus, on passing
to covers, we can ensure that we have thick tubes around a geodesic if
the fundamental group is residually finite.

Again, we have an analogue of this property in the case of topological
geodesics. This result is essentially present in the work of
Gabai~\cite{Ga1}, who basically shows that thick tubes are present
when topological geodesics exist and $\pi_1(M)$ is residually
finite. Thus, the notion of topological geodesics, from the point of
view of thick tubes, is implicit in Gabai's work.

\begin{thm}[See theorem~\ref{T:tubes}] 
Let $\gamma\in\pi_1(M)$ be a primitive element and let $k\in\R$. Then
there is a geodesic $c\subset M$ and a finite cover $M'$ of $M$ such
that $c$ lifts to $M'$ and there is an embedded solid torus $Q$ that
contains $c$ and so that $d(\del Q,c)$ is larger than $k$ and $Q\sm
int(N(c))=T^2\times [0,1]$.
\end{thm} 

\begin{remark} 
Topological geodesics arise naturally while studying elliptic 
$3$-mani-\\folds. For instance, whether homotopy lens spaces are lens
spaces is equivalent to the existence of (the analogue of) a topological
geodesic. Further, lens spaces can be distinguished by considering the
homotopy classes of topological geodesics in a given manifold.
\end{remark}

We now turn to applications. A basic question in topology is to what
extent the homotopy type of a manifold determines the manifold. For
aspherical manifolds, in particular irreducible $3$-manifolds with
infinite fundamental group, conjecturally, pairs of homotopy equivalent
manifolds are always homeomorphic.

One of the fundamental theorems in $3$-manifold topology, due to
Waldhausen, asserts that this is so for so-called Haken
$3$-manifolds. These include irreducible manifolds with non-trivial
boundary. Following Gabai~\cite{Ga1}, we prove rigidity results by
deleting solid tori and reducing to the case of manifolds with
boundary. Partial results along these lines have also been obtained by
Dubois~\cite{Du} who simplified the previous proof by Gabai.

Our next result is that a large class of $3$-manifolds are
\emph{virtually rigid}, i.e., pairs of homotopy equivalent manifolds
have finite covers which are homeomorphic. Gabai~\cite{Ga1} has shown
this assuming residual finiteness and essentially the existence of
topological geodesics.

\begin{theorem}\label{T:rigid}
Suppose $M$ is a closed, irreducible $3$-manifold with $\pi_1(M)$
infinite, residually finite and word-hyperbolic (or
semi-hyperbolic). Then if $f\co M\to N$ is a homotopy equivalence,
there exist finite covers $M'$ and $N'$ of $M$ and $N$ and a lift
$f'\co M'\to N'$ of $f$ which is homotopic to a homeomorphism.
\end{theorem}

\begin{remark} 
We only use word-hyperbolicity to show the tameness of certain covers
of $M$, which also follows under some weaker hypotheses.  
\end{remark}

Our proof (which is already present in Gabai's paper) is a
generalization of that of Gabai asserting the same when $N$ is
hyperbolic. While this has subsequently been strengthened to showing
rigidity for such manifolds (see~\cite{Ga2} and \cite{Ga3}), the
methods used are rather special to hyperbolic manifolds and are
unlikely to generalise.

We also prove, under the same hypothesis, another rigidity result of a
complementary nature. This uses both the existence and uniqueness of
topological geodesics.

\begin{thm}[See theorem~\ref{T:self}]
Let $M$ be irreducible and with word-hyperbolic fundamental group.  If
$f:M\to M$ is a homeomorphism homotopic to the identity then there is
a finite cover $M'$ of $M$ and a lift $f':M'\to M'$ of $f$ such that
$f'$ is isotopic to the identity.
\end{thm}

Further applications of topological geodesics will be studied elsewhere.

\begin{acknowledgements}
We would like to thank Yair Minsky, Darren Long and David Gabai for
helpful comments and conversations.
\end{acknowledgements}

\section{Definition and Existence}

We assume henceforth that $M$ is a closed, irreducible $3$-manifold
with $\pi_1(M)$ word-hyperbolic, residually finite (and infinite). We
shall generalize this to the semi-hyperbolic case, but we begin with
the word-hyperbolic case which is easier. We make (as in the
introduction) the following definition.

\begin{definition} 
An embedded curve $\gamma$ in $M$ is a topological geodesic if a
component $\widetilde\gamma$ of its inverse image in $\widetilde M$ is unknotted.
\end{definition}

The existence of geodesics is based on the following lemma.

\begin{lemma}\label{T:torus}
For $\gamma\in\pi_1(M)$, let $M_\gamma=\widetilde M/\langle g\rangle$ be
the quotient of $\widetilde M$ by the the group of deck transformations
generated by $\gamma$. Then $\widetilde M/\langle g\rangle=S^1\times \R^2$.
\end{lemma}
\begin{proof}
As $\pi_1(M)$ is word-hyperbolic, the universal cover $\widetilde M$ has a
compactification to $B^3$, and the action by deck transformations
extends to $B^3$. The action of $\gamma$ has two fixed points $p$ and
$q$, and $\gamma$ acts properly discontinuously on $B^3\sm\{p,q\}$ with
quotient $D^2 \times S^1$. The result follows as $M_\gamma$ is the
interior of this manifold.
\end{proof}

\begin{theorem}\label{T:exist} 
Let $M$ be an irreducible $3$-manifold with word-hyperbolic\break 
$\pi_1(M)$. Then given $\gamma\in\pi_1(M)$, there is a topological
geodesic $c$ that represents $\gamma$.
\end{theorem}
\begin{proof}
We simply take the image in $M$ of a curve $c$ in $M_\gamma$ that is a
core of the solid torus constructed above. It is easy to ensure that
the image is embedded.
\end{proof}

\begin{remark}  
The above theorem shows that there exists a topological geodesic
representing every element in the fundamental group, rather than each
conjugacy class of elements, which is the case for geodesics in a
Riemannian manifold with negative sectional curvature. But given a
geodesic representing an element, there is an obvious construction of
a geodesic representing any conjugate element. So the real issue is
that the topological geodesic representing a conjugacy class should be
\emph{unique}.
\end{remark}

Next, we weaken the hypothesis on $\pi_1(M)$. 
For the sake of completeness we outline the definition of a semi-hyperbolic group below (see \cite{AB}).
\begin{definition}
For a metric space $(X, d)$ set $P(X)$ for the set of eventually constant maps 
$p:{\Z}_+\to X$, thought as finite discrete paths in $X$. For $p\in P(X)$ 
one denotes by $T_p$ the smallest integer at which $p$ becomes constant.
A {\em bicombing} of $X$ consists   of a choice of a (combing) path $s_{(x,y)}\in P(X)$ 
joining the points $x$ and $y$ of $X$, for all $x,y\in X$. 
\begin{itemize}
\item The bicombing $s$ is {\em quasi-geodesic} if there exist
constants $\lambda, \varepsilon$ such that  $s_{(x,y)}|_{[0, T_{s_{(x,y)}}]}$ 
is a $(\lambda, \varepsilon)$-quasi-geodesic, for all $x,y\in X$. This means that 
\[ d(s_{(x,y)}(t),s_{(x,y)}(t')) \leq \lambda |t-t'| +\varepsilon, \]
is fulfilled for all $x,y\in X$, $t, t'\in {\Z}_+$, thus the combing paths are uniformly closed to geodesics. 
\item The bicombing $s$ is {\em bounded} if there exist constants $k_1\geq 1, k_2\geq 0$
such that for all $x,y,x',y'\in X, t\in {\Z}_+$ one has 
\[ d(s_{(x,y)}(t),s_{(x',y')}(t))\leq k_1 max\{d(x,x'), d(y,y')\} +k_2.\]
\item The metric space $X$, acted upon isometrically by the discrete 
group $\Gamma$, is 
said to be a {\em semi-hyperbolic $\Gamma$-metric space} if it admits a 
bounded quasi-geodesic bicombing $s$ which is $\Gamma$-equivariant i.e.\ 
$\gamma s_{(x,y)}(t)=s_{(\gamma x,\gamma y)}(t)$ holds 
for all $x,y\in X$, $\gamma\in \Gamma$, $t\in {\Z}_+$. 
\item The (finitely generated) group $\Gamma$ is {\em semi-hyperbolic} if $\Gamma$ with the 
word metric (associated to some system of generators) is 
a semi-hyperbolic $\Gamma$-metric space, with respect to the action by left multiplication. 
\end{itemize}
\end{definition}
The semi-hyperbolicity is independent on the system of generators. 
It is known that hyperbolic groups and biautomatic groups 
are semi-hyperbolic.

\begin{lemma}
Assume $\gamma$ is an element of infinite order and $\pi_1(M)$ is a
CAT(0) group (i.e.\ it acts properly discontinuously and co-compactly by
isometries on a geodesic CAT(0) metric space), or more generally 
a semi-hyperbolic group.  
Then the conclusion of
lemma~\ref{T:torus} holds true.
\end{lemma}

\begin{proof}
Any infinite cyclic subgroup ${\Z}\subset \pi_1(M)$ yields a geodesic line in
$X$ by the Flat Torus Theorem (see \cite{AB,V}) and so it is a
quasi-convex subgroup (see \cite{AB}, Thm.9.11 for details).  The
main theorem of \cite{M} implies that $\widetilde
M/{\langle\gamma\rangle}$ is a missing boundary manifold, hence the 
interior of a solid torus. 
In the semi-hyperbolic case an infinite cyclic group lies
in the center of its centralizer (which is a finite extension of the
infinite cyclic group) and the latter is both semi-hyperbolic and a
quasi-convex subgroup of $\pi_1(M)$.
\end{proof}

\begin{corollary}\label{T:exist2}
In the hypotheses of theorem~\ref{T:exist} we may replace
word-hyper-\break bolic by semi-hyperbolic.
\end{corollary}

Thus, we may henceforth weaken our hypothesis on $\pi_1(M)$ to require
only semi-hyperbolicity in place of word hyperbolicity.

Using some results of M.Bridson one can extend the previous 
lemma to nil and sol 3-manifolds, thus holding true for 
all fundamental groups of geometric 3-manifolds. Notice that
M.Kapovich has announced that atoroidal CAT(0) 3-manifold groups are
word hyperbolic, which improves on a previous theorem of L.Mosher.

\section{Thick tubes} 

In this section, we show that a geodesic $c$ corresponding to a
primitive element of $\pi_1(M)$ has a \emph{thick tube} around it in a
finite cover, i.e., a solid torus $Q$ that contains $c$ and so that
$d(\del Q,c)$ is large  and $T\sm int(N(c))=T^2\times [0,1]$. The
precise statement is below. Here, and henceforth, we fix a metric on
$M$ and use the pull-back metric on all its covers

The construction is based on the fact due to Darren Long (see
\cite{L}) that in
word-hyperbolic, or more generally atoroidal, groups, maximal cyclic
groups are separable. 

\begin{quotethm}[D. Long] 
Let $G$ be a residually finite group, $\gamma\in G$ an element that
generates a maximal abelian subgroup, and $S$ a finite set disjoint
from $\langle\gamma\rangle$. Then there exists a subgroup $H\subset G$
with finite index such that $\gamma\in H$ but $S\cap H=\emptyset$.
\end{quotethm}

We construct a thick tube in the cover $\widetilde
M/\langle\gamma\rangle=S^1\times\R^2$. This embeds in a cover where
all other small elements are separated from such a group. In this
cover we have a short curve (geodesic) in one homotopy class, but all
closed curves that are not homotopic to a power of this are long. The
precise statement and proof are below.

\begin{theorem}\label{T:tubes} 
Let $\gamma\in\pi_1(M)$ be a primitive element and let $k\in\R$. Then
there is a geodesic $c\subset M$ and a finite cover $M'$ of $M$ such
that $c$ lifts to $M'$ and there is an embedded solid torus $Q$ that
contains $c$ and so that $d(\del Q,c)> k$ and $Q\sm
int(N(c))=T^2\times [0,1]$.
\end{theorem} 

\begin{proof}
Consider the cover $M_\gamma=\widetilde M/\langle g\rangle$ of $M$ and
pick base-points $p$ of $M$ and $p'$ of $M_\gamma$. As $M_\gamma$ is
homeomorphic to $S^1\times R^2$, there is a curve $c$ and a solid
torus $Q$ embedded in $M_\gamma$, with $c$ based at $p'$, so that $c$
and $Q$ are as required. We shall find a finite cover $M'$ of $M$
which is covered by $M_\gamma$ so that the image of $Q$ embeds in
$M'$.

As $Q$ is compact, there is a finite collection $X$ of inverse images
of $p$ in $M_\gamma$ such that if each of these map to distinct points
in an intermediate cover $M'$ between $M_\gamma$ and $M$, then $Q$
embeds in this cover. Choose paths joining $p'$ to each point in $X$
and let $S'$ be the set of elements in $\pi_1(M)$ that are represented
by the images of these paths. Let $S=S'\cup\{xy^{-1}:x,y\in S\}$.

By the above theorem, there exists a subgroup $H$ separating $\gamma$
from $S$. Let $M'$ be the corresponding cover. It is easy to see that
the elements of $X$ map to different points in $M'$, and hence $Q$
embeds in $M'$.
\end{proof}

\begin{definition} 
A {\em thick tube} around a geodesic $c$ is an embedded solid torus $Q$ that
contains $c$ such that $Q\sm int(N(c))=T^2\times [0,1]$.
\end{definition}

Observe that the components of the inverse image of a thick tube
around a geodesic in $\widetilde M$ are unknotted.

\section{Virtual rigidity}

Gabai's proof of virtual rigidity (see~\cite{Ga1}) now
generalizes. Briefly, if there are knots $K$ and $K'$ in $M$ and $N$
whose complements are irreducible and so that, possibly after changing
$f$ by a homotopy, $f(N(K))\subset N(K')$ and $f(M-int(N(K)))\subset
N-int(N(K'))$, then we can apply Waldhausen's result to the knot
exteriors. Gabai's proof proceeds by showing that given a geodesic
$\gamma$ in $M$ (which plays the role of $K$) with a thick tube around
it, the image of this contains a solid torus $W$ (which plays the role
of $N(K')$), so that $f$ can be deformed to a map that takes
$M-int(N(\gamma))$ to $N-int(W)$. Waldhausen's theorem now shows
rigidity.

As $f$ is a homotopy equivalence between compact manifolds, there
exists $C$ such that $d(x,y)\geq C\implies f(x)\neq f(y)$. By the
methods of the previous section, given $C\in\R$ there is a finite
cover $M'$ of $M$, a geodesic $\gamma\subset M'$ and solid tori
$V_i$,$1\leq i\leq 4$ such that $d(\gamma,\del V_1)>C$ and $d(\del
V_i,\del V_{i+1})>C$, $1\leq i\leq 3$. Let $S_i=\del V_i$. Let
$V_0=\gamma$. Replace $M$ by $M'$ and $N$ by its cover with the same
fundamental group as $M'$ (under the identification
$\pi_1(M)\overset{f_*}=\pi_1(N)$.

The rest of the proof is exactly the same as that of Gabai
(hyperbolicity is not used beyond this stage - only that the tori are
far apart and topologically of a standard form). We outline below the
main steps.

Let $g\co N\to M$ be the homotopy inverse of $f$. Let $K=f(S_2)$ and
$J=N(K)\cup(\text{components of $N\sm K$ disjoint from $f(S_1)\cup
f(S_3)$})$. Then $\del J$ has two components, one of which bounds a
region disjoint from $J$ containing $f(S_1)$ and the other bounds a
region disjoint from $J$ containing $f(S_3)$. Further, $J$ is
irreducible and $[K]$ generates $H_2(J)=\Z$. All this follows from the
fact that $d(x,y)\geq C\implies f(x)\neq f(y)$.
 
Next, $J$ contains a homologically non-trivial torus which bounds in
$N$ a solid torus $W$ containing $f(\gamma)$. This is constructed
using the fact that the Thurston norm equals the singular norm, and
that we have a singular torus $T$. Further, $g\co T\to N\sm int(W)$ and
the inclusion $T\to N\sm int(W)$ induce injections between the
fundamental groups.

Now we can deform $f$ and $g$ so that they restrict to give a homotopy
equivalence between $M\sm V_2$ and $N\sm W$. Waldhausen's theorem now
gives the result.

\section{Uniqueness} 

\begin{lemma}[Engulfing lemma]
Any curve $d$ in $M$ homotopic to a geodesic $c$ is contained in a tube
around the geodesic $c$ in some finite cover.  
\end{lemma}
\begin{proof} 
As $c$ and $d$ are homotopic, there exists an annulus $A$ with
boundary components $c$ and $d$. On passing to a cover with a
sufficiently thick tube $Q$ around $c$, the annulus $A$ is contained
in the tube $Q$. Hence $d$ has a lift that is contained in $Q$.
\end{proof}

\begin{lemma}[Core lemma] 
If $c$ and $c'$ are homotopic geodesics and $c'$ is contained in a
thick tube $Q$ around $c$, then $c'$ is isotopic to $c$.
\end{lemma}
\begin{proof}
Consider $\pi_1(Q\sm c')$. The inverse image of the solid torus $Q$ in
$\widetilde M$ contains an unknotted cylinder $D$ that contains a unique
component $C'$ of the inverse image of $c'$. The group $\pi_1(D\sm
C')=\pi_1(\widetilde M\sm C')$ is the kernel of the map $\pi_1(Q\sm c')\to
\Z$ that maps the longitude of the torus $\del Q$ (hence those of $c$
and $c'$), to $1$ and the meridian to $0$. As $C'$ is an unknot
(because $c'$ is a geodesic) this kernel is $\Z$, hence $\pi_1(Q\sm
c')=\Z^2$. This implies that $c'$ is isotopic to $c$.
\end{proof} 

\begin{theorem}\label{T:unique}
Let $M$ be an irreducible $3$-manifold with $\pi_1(M)$
word-hyper-\break bolic (or semi-hyperbolic) and residually
finite. Suppose $c$ and $c'$ are homotopic topological geodesics in
$M$ representing a primitive class in $\pi_1(M)$. Then there exists a
finite cover $M'$ of $M$ such that $c$ and $c'$ lift to isotopic
curves in $M'$.
\end{theorem}
\begin{proof}
As $c$ and $c'$ are homotopic, there is an annulus $A$ bounding $c$
and $c'$. By passing to a cover $M'$ with a sufficiently thick tube
$Q$ around $c$, we may ensure that $A$ is contained in $Q$. The Core
lemma now implies that $c$ and $c'$ are isotopic.
\end{proof}

\section{Homotopy versus isotopy}
It is known that for an irreducible manifold homotopic self-homeomorphisms 
are isotopic provided that the manifold is 
hyperbolic (a consequence of Gabai's rigidity  \cite{Ga2})
or Seifert fibered (by the Scott theorem) or lens spaces.   
\begin{theorem}\label{T:self}
Let $M$ be irreducible and with word-hyperbolic
fundamental group.  If $f:M\to M$ is a homeomorphism homotopic to the
identity then there is a finite cover $M'$ of $M$ and a lift $f':M'\to
M'$ of $f$ such that $f'$ is isotopic to the identity.
\end{theorem}

\begin{proof} 

Let $\gamma$ be a topological geodesic corresponding to a primitive
class in $\pi_1(M)$. Then the geodesics $f(\gamma)$ and $\gamma$ are
homotopic hence there exists a finite cover $M'$ such that $f(\gamma)$
and $\gamma$ lift to isotopic curves.  Since $f$ is homotopic to
identity there is no obstruction in lifting it to a homotopy
equivalence $f'$ of the finite covering $M'$.  One can further assume
(by means of some isotopy on $M'$) that $f(\gamma')=\gamma'$
(pointwise).

The proof of the core lemma shows that $M'-\gamma'$ is atoroidal hence
hyperbolic (since Haken).  Furthermore the restriction
$f'|_{M'-\gamma'}$ is a homeomorphism.  Therefore $f'|_{M'-\gamma'}$
is homotopic to an isometry of $M'-\gamma'$ (by Mostow rigidity) and
hence isotopic to an isometry (by Waldhausen's theorem).

Let $j$ be this isometry.  Since the isometry group of a hyperbolic
manifold of finite volume is finite it follows that $j$ is of finite
order.  Further $j$ has an extension $g$ (by identity) to all of $M'$,
by asking $g$ to keep pointwise $\gamma'$.  In particular $g$ is a
periodic homeomorphism of $M'$.

Let consider the lift $h$ of $g$ to the universal covering
$\widetilde{M}={\bf R}^3$ (which is also the universal covering of $M'$). 
The action of $h$ extends continuously to the compactification 
(over the boundary sphere) to a homeomorphism of the ball $B^3$. 
The action by deck transformations extends to one by  homeomorphisms of 
the compactification obtained by adding the boundary of the group 
$\pi_1(M)$, because this is word-hyperbolic. In our case the boundary 
is the sphere at infinity. 

But $f$ is homotopic to identity, hence the action induced on the
boundary is trivial.  This shows that $h$ is the identity on the
boundary sphere.

\begin{lemma}
A periodic homeomorphism of $B^3$ which restricts to  identity 
on the boundary is the identity. 
\end{lemma}
\begin{proof}
We will use a theorem of Newman (\cite{Ne}) improved by 
Smith (see \cite{Sm}) in its variant stated in (\cite{Bre},Thm.9.5,
p.157). It states that a compact Lie group acting effectively on a
connected topological manifold has a nowhere dense fixed point set. 

One considers the finite cyclic group action induced by our periodic
homeomorphism on the ball $B^3$. This action extends to the sphere
$S^3$ by the identity outside the upper hemisphere. Then the fixed
point set contains a 3-ball and the previous result shows that the
action cannot be effective, and hence $h$ must be the identity map.
\end{proof} 
\noindent In particular $f'$ is isotopic to the identity.
\end{proof}

\section{Concluding remarks}\label{S:hyp}

While we have assumed word-hyperbolicity (or semi-hyperbolicity) for
the sake of definiteness, we actually need much less. Assuming
$\pi_1(M)$ is atoroidal, we also need everywhere that $\widetilde M=\R^3$,
or equivalently that $\widetilde M^3$ is \emph{tame} (i.e, is the interior
of a compact manifold). Existence and uniqueness use the tameness of
$M_\gamma$. Finally, virtual rigidity needs tameness for \emph{some}
$M_\gamma$.

There are contractible $3$-manifolds, namely Whitehead manifolds,
different from $\R^3$. There are also $3$-manifolds different from
$S^1\times \R^2$ that have fundamental group $\Z$ and universal cover
$\R^3$ (see, for example, \cite{ST}). Conjecturally these manifolds do
not admit free co-compact actions, but so far this is known only under
additional assumptions on the group (for example word-hyperbolicity).

One can generalize the notion of topological geodesics to include
knots having a lift in $\R^3$ whose fundamental group is free (or even
residually nilpotent following \cite{Du}). Although weak versions of
the uniqueness could be proved in this context, the obstruction in
deriving rigidity results is the existence.  The existence of at least
one topological geodesic (homotopically nontrivial) appears to be as
difficult in this more general framework as in the present case.

An even weaker notion of a topological geodesic is one whose lift to
the universal cover is not a satellite knot. In the presence of
topological geodesics, the proof of the Core lemma shows that any
curve that is not a topological geodesic lifts to a finite cover where
it is a satellite knot, and hence is also not a topological geodesic
in this weaker sense.

\bibliographystyle{plain}

\Addresses\recd
\end{document}